\documentclass[12pt,reqno]{amsart}
\usepackage{amssymb,amsthm,amsmath,amstext,amsxtra}
\usepackage{mathrsfs,bm}
\usepackage{fullpage}
\usepackage[all]{xy}
\usepackage{booktabs}
\usepackage{mathtools}
\usepackage{hyperref}
\hypersetup{colorlinks=true,urlcolor=blue,citecolor=blue,linkcolor=blue}
\usepackage{enumerate}
\usepackage{ stmaryrd }
\usepackage{enumitem}
\usepackage{comment}
\usepackage{multirow}
\usepackage{colonequals}
\usepackage{tikz}
\usepackage{marginnote}
\usepackage[export]{adjustbox}

\usepackage{tikz}
\usepackage{tikz-cd}
\usepackage{xcolor}

\numberwithin{equation}{section}

\newtheorem{theorem}{Theorem}[section]

\newtheorem{remark}[theorem]{Remark}

\begin{document}

\title{Rank-2 attractors and Deligne's conjecture}

\author{Wenzhe Yang}
\address{SITP, Physics Department, Stanford University, CA, USA, 94305}
\email{yangwz@stanford.edu}

\begin{abstract}
In this paper, we will study the arithmetic geometry of rank-2 attractors, which are Calabi-Yau threefolds whose Hodge structures admit interesting splits. We will develop methods to analyze the algebraic de Rham cohomologies of rank-2 attractors, and we will illustrate how our methods work by focusing on an example in a recent paper by Candelas, de la Ossa, Elmi and van Straten. We will look at the interesting connections between rank-2 attractors in string theory and Deligne's conjecture on the special values of $L$-functions. We will also formulate several open questions concerning the potential connections between string theory and number theory.
\end{abstract}

%\vspace{-30pt}
\maketitle
\setcounter{tocdepth}{1}
\vspace{-13pt}
\tableofcontents
\vspace{-13pt}

\section{Introduction}
The attractor mechanism occurs in the study of supersymmetric black holes in IIB string theory \cite{Ferrara1995, Ferrara1996}. At low energies, the ten dimensional IIB string theory is well-described by supergravity, which is a classical theory with extended stringy objects approximated by point particles. The attractor mechanism is a very important tool to construct BPS black holes by compactifying the supergravity theory on Calabi-Yau threefolds. In the paper \cite{Moore}, Moore has shown that attractors have many interesting arithmetic properties, which are closely related to algebraic number theory. This is quite surprising, as intuitively the attractors are analytic complex manifolds, hence there is no reason to believe its nature is in fact arithmetic. 

A special kind of attractors is called rank-2 attractor, in which case the middle Hodge structure of the Calabi-Yau threefold admits an interesting split. However, generally it is very difficult to find examples of rank-2 attractors, even by numerical methods. In a recent paper \cite{Candelas}, the authors have found a rank-2 attractor, which does have many interesting arithmetic properties, in particular it is modular. The motivation of this paper is to continue the study of the arithmetic properties of the rank-2 attractors, and find more connections between them and number theory; in particular, their connections with Deligne's conjecture on the special values of $L$-functions. We will also formulate several open questions concerning the possible connections between physics, e.g. string theory and number theory. For example, whether string theory can provide an interpretation to the Deligne's period of the Calabi-Yau threefold. In order to facilitate our discussions, we will need more formal mathematical theories, e.g. the concept of pure motives, a beautiful idea originally from Grothendieck.

\subsection{The rank-2 attractors}

Let us first introduce some properties of the attractor mechanism that will be needed in this paper, while the   readers could consult the paper \cite{Moore} for a more rigorous treatment. We will use the more formal mathematical language, i.e. pure Hodge structure, to explain what is a rank-2 attractor. The readers are referred to the book \cite{PetersSteenbrink} for a quick review of pure Hodge structures. 

In the supergravity theory, the ten dimensional spacetime is taken to be $\mathbb{R}^{3,1} \times X$, where $\mathbb{R}^{3,1}$ is the four dimensional Minkowski spacetime and $X$ is a Calabi-Yau threefold. The charge lattice of the black holes is given by $H^3(X,\mathbb{Z})$ (modulo torsion) \cite{Moore}. In this paper, we will focus on the case where the Hodge number $h^{2,1}$ of $X$ is 1, or equivalently, the dimension of $H^3(X, \mathbb{Q})$ is 4. For simplicity, we will assume $X$ admits a deformation over $\mathbb{P}^1$
\begin{equation} \label{eq:deformationX}
\pi:\mathscr{X} \rightarrow \mathbb{P}^1.
\end{equation}
Let the variable $\varphi$ be a coordinate for the base $\mathbb{P}^1$, then the fiber of the family \ref{eq:deformationX} over a point $\varphi \in \mathbb{P}^1$ will be denoted by $\mathscr{X}_\varphi$. Given a smooth fiber $\mathscr{X}_\varphi$, there is a Hodge decomposition 
\begin{equation} \label{eq:introHodgeDecom}
H^3(X,\mathbb{Q}) \otimes \mathbb{C}=H^{3,0}(\mathscr{X}_\varphi) \oplus H^{2,1}(\mathscr{X}_\varphi) \oplus H^{1,2}(\mathscr{X}_\varphi) \oplus H^{0,3}(\mathscr{X}_\varphi),
\end{equation}
which defines a pure Hodge structure $(H^3(X,\mathbb{Q}),F^p_\varphi)$ on $H^3(X,\mathbb{Q})$ with the Hodge filtration $F^p_\varphi$ given by
\begin{equation}
F^p_\varphi:=\oplus_{q \geq p} H^{q,3-q}(\mathscr{X}_\varphi).
\end{equation}
From now on, the symbol $X$ will also mean the underlying differential manifold of a smooth fiber $\mathscr{X}_\varphi$ in the family \ref{eq:deformationX}. The Hodge filtration $F^p_\varphi$ varies holomorphically with respect to the variable $\varphi$, which in fact defines a holomorphic vector bundle over the smooth locus of the family \ref{eq:deformationX} \cite{KimYang}.  

A point $\varphi \in \mathbb{P}^1$ is called an attractor point if there exists a nonzero charge vector $\gamma \in H^3(X,\mathbb{Z})$ such that the Hodge decomposition of $\gamma$ only has $(3,0)$ and $(0,3)$ components, i.e.
\begin{equation} \label{eq:equationattractorsplit}
\gamma=\gamma^{(3,0)}+\gamma^{(0,3)}.
\end{equation}
A Calabi-Yau threefold $\mathscr{X}_\varphi$ is called an attractor if $\varphi$ is an attractor point. A point $\varphi \in \mathbb{P}^1$ is called a rank-2 attractor point if there exist two linearly independent charge vectors $\gamma_1$ and $\gamma_2$ such that their Hodge decompositions satisfy 
\begin{equation}
\gamma_1=\gamma_1^{(3,0)}+\gamma_1^{(0,3)}, ~\gamma_2=\gamma_2^{(3,0)}+\gamma_2^{(0,3)}.
\end{equation}
A Calabi-Yau threefold $\mathscr{X}_\varphi$ is called a rank-2 attractor if $\varphi$ is a rank-2 attractor point. The rank-2 attractors are defined by rather `analytic' conditions, i.e. a priori there is no reason to believe that $\mathscr{X}_\varphi$ has an algebraic structure. That is why it is very surprising that in the paper \cite{Moore}, Moore conjectures that $\mathscr{X}_\varphi$ should be an algebraic variety defined over a number field.

To proceed, we will need to use some more formal mathematics now! If $\mathscr{X}_\varphi$ is a rank-2 attractor, then the two charge vectors $\gamma_1$ and $\gamma_2$ will define a two dimensional sub-Hodge structure $\mathbf{M}^\textbf{att}_{\text{B}}$, whose Hodge type is $(3,0)+(0,3)$. Here the superscript `\textbf{att}' means the attractor part of $H^3(X,\mathbb{Q})$. Then  the pure Hodge structure $(H^3(X,\mathbb{Q}),F^p_\varphi)$ splits into the following direct sum
\begin{equation} \label{eq:introsplitHodge}
(H^3(X,\mathbb{Q}),F^p_\varphi)=\mathbf{M}^\textbf{att}_{\text{B}} \oplus \mathbf{M}^\textbf{ell}_{\text{B}},
\end{equation}
where $\mathbf{M}^\textbf{ell}_{\text{B}}$ with Hodge type $(2,1)+(1,2)$. Here the superscript `\textbf{ell}' means the elliptic part of  $H^3(X,\mathbb{Q})$. As the name has suggested, $\mathbf{M}^\textbf{ell}_{\text{B}}$ is closely related to the theory of elliptic curves.

\subsection{The pure motives and Deligne's conjecture}

In order to study the arithmetic properties of rank-2 attractors, we will need to introduce the more abstract concept of pure motives. Pure motives are very foundational objects in the study of arithmetic algebraic geometry. This concept certainly sounds very difficult at the beginning, but even physicists will appreciate its simplicity and beauty after becoming familiar with it. More specifically, it will greatly facilitate our discussions about the arithmetic properties of rank-2 attractors.

In this paper, we will use the pure motive as a black box, and focus on its classical realizations. The interested readers are referred to the papers \cite{Shamit, KimYang, Yang1} for more formal treatments. Given a smooth algebraic variety $Y$ defined over $\mathbb{Q}$, i.e. $Y$ is defined by polynomial equations with rational coefficients, then a pure motive 
\begin{equation}
\mathbf{M}:=h^n(Y),~ n \in \mathbb{Z}
\end{equation}
has three important classical realizations:
\begin{enumerate}
\item The Betti realization $\mathbf{M}_\text{B}$. Notice that the variety $Y$ has a complex manifold structure. The Betti realization $\mathbf{M}_\text{B}$ is just the singular cohomology group $H^n(Y,\mathbb{Q})$. Moreover, there is a pure Hodge structure on $H^n(Y,\mathbb{Q})$ induced by the Hodge decomposition 
\begin{equation}
H^n(Y,\mathbb{Q}) \otimes \mathbb{C}=H^{n,0}(Y) \oplus H^{n-1,1}(Y)  \oplus \cdots \oplus H^{1,n-1}(Y) \oplus H^{0,n}(Y),
\end{equation}
where the Hodge filtration $F^pH^n(Y,\mathbb{C})$ is given by
\begin{equation} \label{eq:HodgefiltrationHnY}
F^pH^n(Y,\mathbb{C}) =\oplus_{q \geq p} H^{q,n-q}(Y).
\end{equation}
Together with this pure Hodge structure, the Betti realization gives us the Hodge realization of $\mathbf{M}$.

\item The de Rham realization $\mathbf{M}_\text{dR}$, which is the algebraic de Rham cohomology group $H^n_{\text{dR}}(Y)$ defined by the algebraic forms on $Y$ \cite{KimYang, Nekovar}. The group $H^n_{\text{dR}}(Y)$ has a Hodge filtration $F^p(H^n_{\text{dR}}(Y))$, the construction of which can be found in the papers \cite{Nekovar, Yang1}. There is a canonical comparison isomorphism between the Betti realization and the de Rham realization, which is induced by the integration of the algebraic forms over the homological cycles. Under this comparison isomorphism, the Hodge filtration $F^p(H^n_{\text{dR}}(Y))$ corresponds to the filtration $F^pH^n(Y,\mathbb{Q})$ \cite{Nekovar}.

\item The \'etale realization $\mathbf{M}_\ell$, which is the \'etale cohomology group $H^n_{\text{\'et}}(Y,\mathbb{Q}_\ell)$ of $Y$, where $\ell$ is a prime number. It is a continuous representation of the absolute Galois group $\text{Gal}(\overline{\mathbb{Q}}/\mathbb{Q})$ \cite{DeligneL, Nekovar, Yang1}. The \'etale realization is crucial for the study of the arithmetic properties of $Y$, e.g. zeta functions and modularity. The good news is that in this paper, we will only use some properties of $\mathbf{M}_\ell$, and we do not need to dig deep into the technical details.
\end{enumerate}
It is certainly very safe and also very helpful to comprehend the concept of pure motives through their classical realizations, which is also the point of view adopted in this paper.

Let us now look at an example \cite{Shamit, KimYang, Yang1}. The Tate motive $\mathbb{Q}(1)$ is by definition the dual of the Lefschetz motive $h^2(\mathbb{P}^1_{\mathbb{Q}})$, whose classical realizations are:
\begin{enumerate}
\item $\mathbb{Q}(1)_{\text{B}}=(2 \pi i)\, \mathbb{Q}$, which admits a pure Hodge structure of the type $(-1,-1)$.

\item $\mathbb{Q}(1)_{\text{dR}}=\mathbb{Q}$, with the Hodge filtration given by $F^0=0$ and $F^{-1}=\mathbb{Q}$.

\item $\mathbb{Q}(1)_{\ell}= \varprojlim_n \mu_{\ell^n}(\overline{\mathbb{Q}}) \otimes_{\mathbb{Z}_{\ell}} \mathbb{Q}_\ell$. Here $\mu_{\ell^n}(\overline{\mathbb{Q}})$ consists of the $\ell^n$-th root of unity which admits a natural action by $\mathbb{Z}/\ell \mathbb{Z}$.
\end{enumerate}
The Tate motive $\mathbb{Q}(n)$ is the $n$-fold tensor product $\mathbb{Q}(1)^{\otimes n}$ \cite{KimYang, Nekovar}.  The Tate twist of $\mathbf{M}$ by $\mathbb{Q}(n)$, $n \in \mathbb{Z}$, i.e. $\mathbf{M} \otimes \mathbb{Q}(n)$, will be denoted by $\mathbf{M}(n)$. Here the tensor products of pure motives behave similarly as the tensor products of classical cohomology theories.

The \'etale realization $\mathbf{M}_\ell$ of the pure motive $\mathbf{M}$ is a continuous representation of the absolute Galois group $\text{Gal}(\overline{\mathbb{Q}}/\mathbb{Q})$, which allows us to define an $L$-function $L(\mathbf{M},s)$ for $\mathbf{M}$. While the details of the construction of $L(\mathbf{M},s)$ are left to the papers \cite{Nekovar, Yang1}. The $L$-function $L(\mathbf{M},s)$ satisfies the following property
\begin{equation}
L(\mathbf{M}(n),s)=L(\mathbf{M},n+s).
\end{equation}
Therefore, in order to study the special value of $L(\mathbf{M},s)$ at $s=n$, it is necessary and sufficient to study the special value of $L(\mathbf{M}(n),s)$ at $s=0$. The pure motive $\mathbf{M}$ is called critical if and only if its Hodge numbers satisfy certain conditions \cite{DeligneL}. For example, if the Hodge type of $\mathbf{M}$ is $(3, 0)+(0, 3)$, then $\mathbf{M}(n)$ is critical if and only if $n=1, 2$. For the pure motive $\mathbf{M}(n)$, the Deligne's period $c^+(\mathbf{M}(n))$ can be computed from its Betti realization $\mathbf{M}_\text{B}$ and de Rham realization  $\mathbf{M}_\text{dR}$. Now, we are ready to state a deep  conjecture of Deligne.

\vspace*{0.1in}

\textbf{Deligne's conjecture}: \textit{If the pure motive $M(n)$ is critical, then the value $L(M(n),0)$ is a rational multiple of the period $c^+(M(n))$.}

\vspace*{0.1in}

\subsection{The arithmetic properties of rank-2 attractors}

Let us now look at the arithmetic properties of rank-2 attractors. In the paper \cite{Moore}, Moore has conjectured that rank-2 attractors are algebraically defined over number fields. But for simplicity, in this paper we will focus on the rank-2 attractors that are algebraic varieties defined over $\mathbb{Q}$. More precisely, we will focus on the case where the deformation of $X$ in \ref{eq:deformationX} is algebraically defined over $\mathbb{Q}$, and we will only consider the rational rank-2 attractor points $\varphi$, i.e. $\varphi \in \mathbb{Q}$. Then the Hodge conjecture tells us that over a number field $K$, the pure motive $h^3(\mathscr{X}_\varphi)$ of a rank-2 attractor $\mathscr{X}_\varphi, \varphi \in \mathbb{Q}$ splits into the following direct sum 
\begin{equation} \label{eq:Hodgeconjecturespli}
h^3(\mathscr{X}_\varphi)=\mathbf{M}^\textbf{att} \oplus \mathbf{M}^\textbf{ell},
\end{equation}
where the Hodge realization of the sub-motive $\mathbf{M}^\textbf{att}$ (resp. $\mathbf{M}^\textbf{ell}$) is the Hodge structure $\mathbf{M}^\textbf{att}_{\text{B}} $ (resp. $\mathbf{M}^\textbf{ell}_{\text{B}} $) in the formula \ref{eq:introsplitHodge} \cite{Shamit}. Furthermore, if we assume $K$ is the rational field $\mathbb{Q}$, then the split \ref{eq:Hodgeconjecturespli} implies that the four dimensional Galois representation $H^3_{\text{\'et}}(\mathscr{X}_\varphi,\mathbb{Q}_\ell)$ splits into the direct sum
\begin{equation}
H^3_{\text{\'et}}(\mathscr{X}_\varphi,\mathbb{Q}_\ell)=\mathbf{M}^\textbf{att}_\text{\'et} \oplus \mathbf{M}^\textbf{ell}_\text{\'et},
\end{equation}
where $\mathbf{M}^\textbf{att}_\text{\'et}$ (resp. $\mathbf{M}^\textbf{ell}_\text{\'et}$) is the \'etale realization of $\mathbf{M}^\textbf{att} $ (resp. $\mathbf{M}^\textbf{ell} $).

The Tate twist $\mathbf{M}^\textbf{ell} \otimes \mathbb{Q}(1)$ is of Hodge type $(1,0)+(0,1)$ \cite{Shamit, KimYang, Yang1}. From the modularity theorem of elliptic curves, we deduce that the Galois representation $\mathbf{M}^\textbf{ell}_\text{\'et} \otimes \mathbb{Q}_\ell(1)$ is modular \cite{Diamond}. More concretely, at a good prime number $p$, the characteristic polynomial of the geometric Frobenius for $\mathbf{M}^\textbf{ell}_\text{\'et} \otimes \mathbb{Q}_\ell(1)$ is of the form
\begin{equation}
1-a_pT+pT^2,
\end{equation}
where $a_p$ is the $p$-th coefficient of the $q$-expansion of a weight-2 newform \cite{Diamond}. While for $\mathbf{M}^\textbf{ell}_\text{\'et} $, the characteristic polynomial of the geometric Frobenius at $p$ will be
\begin{equation}
1-a_p(pT)+p(pT)^2.
\end{equation}
From the paper \cite{Yui}, $\mathbf{M}^\textbf{att}_\text{\'et}$ is also modular, i.e. at a good prime number $p$, the characteristic polynomial of the geometric Frobenius is of the form
\begin{equation}
1-b_pT+p^3T^2,
\end{equation}
where $b_p$ is the $p$-th coefficient of the $q$-expansion of a weight-4 Hecke eigenform. On the other hand, the pure motives $\mathbf{M}^\textbf{att} \otimes \mathbb{Q}(2)$ and $\mathbf{M}^\textbf{att} \otimes \mathbb{Q}(1)$ are both critical, therefore a very interesting question is to check whether they satisfy Deligne's conjecture or not \cite{DeligneL,Yang1}.

In the paper \cite{Candelas}, the authors have found a rank-2 attractor defined over $\mathbb{Q}$, which is a smooth fiber over the rational point $\varphi=-1/7$ in a one-parameter family of Calabi-Yau threefolds. This rank-2 attractor will be denoted by $\mathscr{X}_{-1/7}$. The authors of \cite{Candelas} have numerically computed the zeta functions of $\mathscr{X}_{-1/7}$ for small prime numbers. They find that at a good prime number $p$, the zeta function of the pure motive $h^3(\mathscr{X}_{-1/7})$ is of the form
\begin{equation}
(1-a_p(pT)+p(pT)^2)(1-b_pT+p^3T^2).
\end{equation}
Here $a_p$ is the $p$-th coefficient of the $q$-expansion of a weight-2 modular form $f_2$ for the modular group $\Gamma_0(14)$, which is labeled as \textbf{14.2.a.a} in LMFDB. While $b_p$ is the $p$-th coefficient of the $q$-expansion of a weight-4 modular form $f_4$ also for the modular group $\Gamma_0(14)$, which is labeled as \textbf{14.4.a.a} in LMFDB \cite{Candelas}. The authors of \cite{Candelas} have also numerically computed the special values of the $L$-functions, $L(f_2,1)$, $L(f_4,1)$ and $L(f_4,2)$. Besides, they find a special number $v^\perp$ such that the $j$-value of $\frac{1}{2}+i \,v^\perp$ is $215^3/28^3$, which corresponds to the modular curve $X_0(14)$ associated to the modular group $\Gamma_0(14)$. The weight-2 newform associated to $X_0(14)$ under the modularity theorem of elliptic curves is just the modular form $f_2$.

In the paper \cite{Yang1}, a new method to compute the Deligne's periods for Calabi-Yau threefolds in a mirror family has been developed. Based on this method, it has been numerically shown that the critical motive $h^3(\mathscr{X}_{-1/7} \otimes \mathbb{Q}(2)) $ does satisfy the prediction of Deligne's conjecture. This paper can be considered as a further development of the paper \cite{Yang1}. In this paper, we will look at the algebraic de Rham cohomology of rank-2 attractors, which sheds further light on the split \ref{eq:Hodgeconjecturespli}. More precisely, we will explicitly construct the de Rham realization of the split \ref{eq:Hodgeconjecturespli} for the rank-2 attractor $\mathscr{X}_{-1/7}$. Then we will compute the Deligne's periods for the critical motives $\mathbf{M}^\textbf{att} \otimes \mathbb{Q}(2)$, $\mathbf{M}^\textbf{att} \otimes \mathbb{Q}(1)$ and $\mathbf{M}^\textbf{ell} \otimes \mathbb{Q}(2)$, and we will numerically show they do satisfy Deligne's conjecture. Furthermore, we will construct a Betti realization of $\mathbf{M}^\textbf{ell}_{\text{B}}$, which is a two dimensional subspace of $H^3(X,\mathbb{Q})$. We have also found a rank-2 sub-lattice of $H^3(X,\mathbb{Z})$, which induces a lattice structure for $\mathbf{M}^\textbf{ell}_{\text{B}}$. More concretely, we have found a rank-2 sub-lattice of $H^3(X,\mathbb{Z})$, the complex elliptic curve associated to which has the same $j$-invariant as $X_0(14)$.

The outline of this paper is as follows. In Section \ref{sec:mirrorsymmetryattractor}, we will give an overview of mirror symmetry and introduce the attractor equation. In Section \ref{sec:exampleAESZ34}, we will review some results of the paper \cite{Candelas} concerning the rank-2 attractor $\mathscr{X}_{-1/7}$. In Section \ref{sec:deRhamsplit}, we will study the algebraic de Rham cohomology of the rank-2 attractor $\mathscr{X}_{-1/7}$, and show it also admits a split. In Section \ref{sec:verificationDeligne}, we will compute the Deligne's periods for the critical motives  $\mathbf{M}^\textbf{att} \otimes \mathbb{Q}(2)$, $\mathbf{M}^\textbf{att} \otimes \mathbb{Q}(1)$ and $\mathbf{M}^\textbf{ell} \otimes \mathbb{Q}(2)$, and we will numerically verify that they satisfy the predictions of Deligne's conjecture. In Section \ref{sec:constructionEllipticsubmotive}, we will find a rank-2 sub-lattice of $H^3(X,\mathbb{Z})$, and compute the $j$-invariant of the complex elliptic curve associated to it. In Section \ref{sec:conclusion}, we will discuss the conclusions of this paper and some further open questions, which can further our understanding of the arithmetic properties of rank-2 attractors.

\section{Mirror symmetry and the attractor equation} \label{sec:mirrorsymmetryattractor}

In this section, we will briefly review some results of the mirror symmetry of Calabi-Yau threefolds  \cite{Yang1}. We will focus on the one-parameter mirror pairs of Calabi-Yau threefolds, while the details are left to the papers \cite{PhilipXenia, CoxKatz, MarkGross, KimYang}. Given a mirror pair $(X^\vee,X)$ of Calabi-Yau threefolds, here one-parameter means that their Hodge numbers satisfy
\begin{equation}
h^{1,1}(X^\vee)=h^{2,1}(X)=1.
\end{equation}

\subsection{Picard-Fuchs equation}

In this paper, we will further assume that the mirror threefold $X$ admits an algebraic deformation defined over $\mathbb{Q}$ of the form \cite{KimYang}
\begin{equation} \label{eq:mirrorfamily}
\pi:\mathscr{X} \rightarrow \mathbb{P}^1_{\mathbb{Q}},
\end{equation}
where the coordinate of the base variety $\mathbb{P}^1_{\mathbb{Q}}$ has been chosen to be $\varphi$. The fiber of this family \ref{eq:mirrorfamily} over a point $\varphi \in \mathbb{P}^1_{\mathbb{Q}}$ is denoted by $\mathscr{X}_\varphi$. From now on, $X$ will also mean the underlying differential manifold of a smooth fiber of this family. We will also assume that for every smooth fiber $\mathscr{X}_\varphi$, there exists a nowhere-vanishing algebraic threeform $\Omega_\varphi$, which varies holomorphically with respect to $\varphi$ \cite{KimYang}. Moreover, as a threeform on $\mathscr{X}$, $\Omega$ is defined over $\mathbb{Q}$; in particular, if $\varphi$ is rational, then $\Omega_\varphi$ is also defined over $\mathbb{Q}$ \cite{CoxKatz, MarkGross, KimYang}. From the Griffiths transversality, $\Omega_\varphi$ satisfies a fourth-order Picard-Fuchs equation
\begin{equation} \label{eq:PFequation}
\mathscr{L}\,\Omega_\varphi=0,
\end{equation}
where $\mathscr{L}$ is a differential operator with polynomial coefficients $R_i(\varphi) \in \mathbb{Q}[\varphi]$
\begin{equation} \label{eq:PicardFuchsOperator}
\mathscr{L}=R_4(\varphi) \,\vartheta^4+R_3(\varphi)\, \vartheta^3 +R_2(\varphi)\, \vartheta^2+ R_1(\varphi) \, \vartheta^1+R_0(\varphi), ~\text{with}~ \vartheta=\varphi \,\frac{d}{d\varphi}.
\end{equation}

We will also assume the point $\varphi=0$ is the large complex structure limit, which means that the monodromy of $\mathcal{L}$ at it is maximally unipotent. More concretely, there exists a small disc $\Delta$ of $\varphi=0$ such that the Picard-Fuchs equation \ref{eq:PFequation} has four canonical solutions of the form 
\begin{equation} \label{eq:PeriodsCan}
\begin{aligned}
\varpi_0 &= f_0,  \\
\varpi_1 &=\frac{1}{2\pi i}\left(f_0 \log \varphi+f_1\right), \\
\varpi_2 &=\frac{1}{(2\pi i)^2}\left( f_0 \log^2 \varphi +2\, f_1 \log \varphi + f_2\right), \\
\varpi_3 &=\frac{1}{(2 \pi i)^3} \left( f_0 \log^3 \varphi +3 \, f_1 \log^2 \varphi +3\, f_2 \log \varphi +f_3 \right),
\end{aligned} 
\end{equation}
where $\{f_j\}_{j=0}^3$ are power series in $\mathbb{Q}[[\varphi]]$ that converge on $\Delta$. If we impose the following condition on $f_j$
\begin{equation} \label{eq:boundarycondition}
f_0(0)=1,~f_1(0)=f_2(0)=f_3(0)=0,
\end{equation}
then the four canonical solutions \ref{eq:PeriodsCan} are unique \cite{KimYang, Yang}. The canonical period vector $\varpi$ is the column vector defined by
\begin{equation}
\varpi:=(\varpi_0,\,\varpi_1,\,\varpi_2,\,\varpi_3)^\top.
\end{equation}
\begin{remark}
In this paper, the multi-valued holomorphic function $\log \varphi$ is chosen to satisfy
\begin{equation}
\log (1)=0,~\log (-1)=\pi i.
\end{equation}
\end{remark}

The algebraic de Rham cohomology $H^3_{\text{dR}}(\mathscr{X}_\varphi)$ of a smooth rational fiber $\mathscr{X}_\varphi,\varphi \in \mathbb{Q}$ is completely determined by the threeform $\Omega_\varphi$ and its derivatives:
\begin{equation} \label{eq:HodgefiltrationMS}
\begin{aligned}
F^3(H^3_{\text{dR}}(\mathscr{X}_\varphi))&=\langle \Omega_\varphi \rangle, \\
F^2(H^3_{\text{dR}}(\mathscr{X}_\varphi))&=\langle \Omega_\varphi,\Omega'_\varphi \rangle, \\
F^1(H^3_{\text{dR}}(\mathscr{X}_\varphi))&=\langle \Omega_\varphi,\Omega'_\varphi,\Omega''_\varphi \rangle, \\
F^0(H^3_{\text{dR}}(\mathscr{X}_\varphi))&=\langle \Omega_\varphi,\Omega'_\varphi,\Omega''_\varphi ,\Omega'''_\varphi \rangle;
\end{aligned}
\end{equation}
where we have used Griffiths transversality \cite{KimYang, Yang1}. Here $\langle \Omega_\varphi \rangle$ means the rational vector space spanned by $\Omega_\varphi$; and $\Omega'_\varphi$ means the first derivative of $\Omega_\varphi$ with respect to $\varphi$, etc. The readers are referred to the papers \cite{KimYang, Yang, Yang1} for more details.

\subsection{Mirror symmetry}

From the Poincar\'e duality, there exists a unimodular skew symmetric pairing on $H_3(X,\mathbb{Z})$ (modulo torsion), which allows us to choose an integral symplectic basis $\{A_0,A_1,B_0,B_1\}$ that satisfy the following intersection pairing \cite{PhilipXenia,CoxKatz,MarkGross}
\begin{equation}
A_a \cdot A_b=0,~~B_a \cdot B_b=0,~~A_a \cdot B_b= \delta_{ab}.
\end{equation}
Suppose the dual of this basis is denoted by $\{\alpha^0,\alpha^1,\beta^0,\beta^1\}$, i.e. we have
\begin{equation}
\alpha^a (A_b)=\delta_{ab}, ~ \beta^a(B_b)=\,\delta_{ab}, \alpha^a(B_b)=\beta^a(A_b)=0,
\end{equation}
then they form a basis for $H^3(X,\mathbb{Z})$ (modulo torsion). From the Poincar\'e duality, we have
\begin{equation} \label{eq:poincarepairing}
\int_X \alpha^a \smile \beta^b=\delta_{ab},~\int_X \alpha^a \smile \alpha^b=0,~\int_X \beta^a \smile \beta^b=0,
\end{equation}
where $\alpha^a \smile \beta^b$ means the cup product between the cohomological  cycles $\alpha^a$ and $\beta^b$, etc \cite{Hatcher}.

\begin{remark}
The torsion of the homology or cohomology groups are irrelevant to this paper, hence they will be ignored.
\end{remark}

The integral periods of the mirror family \ref{eq:mirrorfamily} are given by the integration of the threeform $\Omega_\varphi$ over the symplectic basis $\{A_a,B_a\}_{a=0}^1$ of $H_3(X,\mathbb{Z})$
\begin{equation} \label{eq:IntegralPeriodDefinition}
z_a(\varphi)= \int_{A_a} \Omega_\varphi,~\mathcal{G}_b(\varphi)=\int_{B_b} \Omega_\varphi,
\end{equation}
which are multi-valued holomorphic functions \cite{PhilipXenia,CoxKatz,MarkGross}. The integral period vector $\amalg(\varphi)$ is a column vector defined by
\begin{equation}
\amalg(\varphi)=(\mathcal{G}_0(\varphi),\mathcal{G}_1(\varphi),z_0(\varphi),z_1(\varphi))^\top.
\end{equation}
The integral period vector $\amalg$ forms another basis for the solution space of the Picard-Fuchs equation \ref{eq:PFequation}, hence there exists a transformation matrix $S \in \text{GL}(4,\mathbb{C})$ such that
\begin{equation} \label{eq:PiSomega}
\amalg=S\cdot \varpi.
\end{equation}
The transformation matrix $S$ will be crucial to this paper, and it can be determined by mirror symmetry \cite{KimYang, Yang1}. For later convenience, let us also define the row vector $\beta$ by
\begin{equation} \label{eq:basisbeta}
\beta:=(\beta^0,\beta^1,\alpha^0,\alpha^1),
\end{equation}
which forms a basis for $H^3(X,\mathbb{Z})$. Under the comparison isomorphism between the Betti and the de Rham cohomology, $\Omega_\varphi$ has an expansion of the form
\begin{equation} \label{eq:expansionOmega}
\Omega_\varphi=\beta \cdot \amalg(\varphi)=\mathcal{G}_0(\varphi) \, \beta^0 +\mathcal{G}_1(\varphi) \,\beta^1+z_0(\varphi)\, \alpha^0+z_1(\varphi)\, \alpha^1.
\end{equation}

In all examples of one-parameter mirror pairs, there always exists an integral symplectic basis $\{A_a,B_a\}_{a=0}^1$ of $H_3(X,\mathbb{Z})$ such that 
\begin{equation} \label{eq:zivarpii01}
z_j(\varphi)=\lambda (2 \pi i)^3 \,\varpi_j(\varphi),~j=0,1; \lambda \in \mathbb{Q}^\times,
\end{equation}
where $\lambda$ is a nonzero rational number \cite{PhilipXenia, CoxKatz, MarkGross}. The mirror map is given by the quotient
\begin{equation}
t=\frac{z_1}{z_0}=\frac{\varpi_1}{\varpi_0}=\frac{1}{2 \pi i}\,\log \varphi+\frac{f_1(\varphi)}{f_0(\varphi)}.
\end{equation}
Near the large complex structure limit $\varphi=0$, formula \ref{eq:boundarycondition} immediately that implies 
\begin{equation}
t=\frac{1}{2 \pi i} \,\log \varphi+ \mathcal{O}(\varphi),
\end{equation}
therefore the large complex structure limit corresponds to $t=i\infty $. In mirror symmetry, the prepotential $\mathcal{F}$ admits an expansion of the form \cite{PhilipXenia,CoxKatz,KimYang}
\begin{equation} \label{eq:Prepotential}
\mathcal{F}=-\frac{1}{6}\, Y_{111}\, t^3 -\frac{1}{2}\, Y_{011}\,t^2-\frac{1}{2}\,Y_{001}\, t-\frac{1}{6}\,Y_{000}+\mathcal{F}^{\text{np}},
\end{equation}
where $\mathcal{F}^{\text{np}}$ is the non-perturbative instanton correction 
\begin{equation}
\mathcal{F}^{\text{np}}=\sum_{n=1}^{\infty} a_n \exp 2 \pi i \,nt.
\end{equation}
The coefficient $Y_{111}$ in the expansion \ref{eq:Prepotential} is the topological intersection number \cite{PhilipXenia,CoxKatz, MarkGross}
\begin{equation}
Y_{111}=\int _{X^\vee} e\wedge e \wedge e,
\end{equation}
where $e$ forms a basis for $H^2(X^\vee,\mathbb{Z})$. Moreover, we will require $e$ to lie in the K\"ahler cone of $X^\vee$, hence $Y_{111}$ will be a positive integer \cite{CoxKatz, MarkGross}. The coefficients $Y_{011}$ and $Y_{001}$ are rational numbers, however the computations of them are more tricky \cite{CoxKatz,KimYang}. In all examples of mirror pairs, the coefficient $Y_{000}$ is always of the form \cite{PhilipXenia}
\begin{equation} \label{eq:PhysicistsY000}
Y_{000}=-3\, \chi (X^\vee)\, \frac{\zeta(3)}{(2 \pi i)^3},
\end{equation}
where $\chi(X^\vee)$ is the Euler characteristic of $X^\vee$. A study of the occurrence of $\zeta(3)$ from the motivic point of view is presented in the paper \cite{KimYang}. Using mirror symmetry, it can be shown that the transformation matrix $S$ in the formula \ref{eq:PiSomega} is of the form \cite{KimYang, Yang1}
\begin{equation} \label{eq:smatrix}
S\,=\lambda(2 \pi i)^3 \,
\left(
\begin{array}{cccc}
 -\frac{1}{3}\, Y_{000} & -\frac{1}{2} \,Y_{001} & 0 & \frac{1}{6}\, Y_{111} \\
 -\frac{1}{2} \,Y_{001} & -\,Y_{011} & -\frac{1}{2} \,Y_{111} & 0 \\
 1 & 0 & 0 & 0 \\
 0 & 1 & 0 & 0 \\
\end{array}
\right),~\lambda \in \mathbb{Q}^\times.
\end{equation}

\subsection{The attractor equation}

Now we are ready to write down the attractor equation in IIB string theory for a nonzero black-hole charge $\gamma \in H^3(X,\mathbb{Z})$ \cite{Moore}. The readers could consult the paper \cite{Moore} for a thorough review of the attractor mechanism and its physical significance. Given a point $\varphi\in \mathbb{P}^1_\mathbb{Q}$, as the dimension of $H^{3,0}(\mathscr{X}_\varphi)$ is $1$, we immediately deduce that the component $\gamma^{3,0}$ in the Hodge decomposition of $\gamma$ satisfies
\begin{equation}
\gamma^{3,0}=C\,\Omega_\varphi=C\,\beta \cdot \Pi(\varphi)=C\,\beta\cdot S \cdot \varpi(\varphi),
\end{equation}
where $C \in \mathbb{C}$ is a nonzero constant. Take its complex conjugation, we obtain
\begin{equation}
\gamma^{0,3}=\overline{\gamma^{3,0}}=\overline{C}\,\beta \cdot \overline{\Pi}(\varphi)=\overline{C}\,\beta\cdot \overline{S} \cdot \overline{\varpi}(\varphi).
\end{equation}
Suppose the expansion of $\gamma$ with respect to the basis $\beta$ is
\begin{equation}
\gamma=\beta \cdot P,\,P\in \mathbb{Z}^4-(0,0,0,0).
\end{equation}
If $\gamma$ only has $(3,0)$ and $(0,3)$ components in its Hodge decomposition (at the point $\varphi$), then we must have
\begin{equation}
P=C\, S \cdot \varpi(\varphi)+\overline{C}\, \overline{S} \cdot \overline{\varpi}(\varphi).
\end{equation}
This is called the attractor equation. The point $\varphi$ is an attractor point if and only if there exists a non-zero constant $C \in \mathbb{C}$ such that 
\begin{equation}
C\, S \cdot \varpi(\varphi)+\overline{C}\, \overline{S} \cdot \overline{\varpi}(\varphi) \in \mathbb{Z}^4-(0,0,0,0).
\end{equation}
Similarly, $\varphi$ is a rank-2 attractor point if and only if there exist two nonzero constants $C_1$ and $C_2$ such that $C_1\, S \cdot \varpi(\varphi)+\overline{C}_1\, \overline{S} \cdot \overline{\varpi}(\varphi) $ and $C_2\, S \cdot \varpi(\varphi)+\overline{C}_2\, \overline{S} \cdot \overline{\varpi}(\varphi) $ are two linearly independent elements in $\mathbb{Z}^4$. 

\begin{remark}
At an arbitrary point $\varphi \in \mathbb{P}^1_\mathbb{Q}$, the numerical value of $\varpi_i(\varphi)$ can be evaluated by numerically solving the Picard-Fuchs equation \ref{eq:PFequation}. The readers can consult the paper \cite{Yang} for more details about this method. However, it is still very challenging to search for rank-2 attractor points numerically using softwares like Mathematica.
\end{remark}

\section{An example of the rank-2 attractor} \label{sec:exampleAESZ34}

In this section, we will review the results of the paper \cite{Candelas} that will be needed in later sections. In  \cite{Candelas}, the authors have found two rank-2 attractors that are defined over $\mathbb{Q}$. The two rank-2 attractors are very similar, e.g. they have the same zeta functions. Therefore we will focus on one example, while the analysis of the other is exactly the same, hence will not be repeated.

\subsection{Attractor equations} \label{sec:attractorequation}

In the paper \cite{Candelas}, the authors have constructed a one-parameter mirror pair $(X^\vee,X)$ of Calabi-Yau threefolds, where the Hodge diamond of the mirror threefold $X$ is of the form 
\begin{center} 
\begin{tabular}{ c c c c c c c }
 &  &  & 1 &  &  &  \\ 
 &  & 0&   & 0&  &  \\   
 & 0&  & 9 &  & 0&   \\  
1&  & 1 &  & 1 & & 1 \\ 
 & 0&  & 9 &  & 0&   \\ 
 &  & 0&   & 0&  &  \\   
 &  &  & 1 &  &  &  \\
\end{tabular}.
\end{center}
There exists an algebraic deformation of $X$ defined over $\mathbb{Q}$
\begin{equation} \label{eq:aeszfamilies}
\pi:\mathscr{X} \rightarrow \mathbb{P}^1_{\mathbb{Q}}.
\end{equation}
Notice that $X$ will also mean the underlying differential manifold of a smooth fiber of this family. The construction of $X$ and its deformation will not be needed in this paper, hence are left to \cite{Candelas}. There exists a family of threefroms $\Omega_\varphi$ for the deformation \ref{eq:aeszfamilies} that satisfies all the assumptions in Section \ref{sec:mirrorsymmetryattractor}. The Picard-Fuchs equation satisfied by $\Omega_\varphi$ is $\mathcal{D}\,\Omega_\varphi=0$ with
\begin{equation} \label{eq:aesz34PF}
\begin{aligned}
\mathcal{D}=& \theta^4 -\varphi(35 \theta^4+70 \theta^3+63 \theta^2+28 \theta +5) + \varphi^2 (\theta+1)^2(259 \theta^2+518 \theta +285)\\
                       &-225 \varphi^3 (\theta+1)^2(\theta+2)^2,~ \theta=\varphi \frac{d}{d \varphi}.
\end{aligned}
\end{equation}
The Picard-Fuchs operator $\mathcal{D}$ has five regular singularities at the points 
\begin{equation}
\varphi=0, 1/25, 1/9, 1, \infty,
\end{equation}
while $\varphi=0$ is the large complex structure limit. The canonical period $\varpi_0$ is given by
\begin{equation}
\varpi_0=1+\sum_{n=1}^\infty a_n \varphi^n;~~a_n=\sum_{i+j+k+l+m=n}\left(\frac{n! }{i!j!k!l!m!} \right)^2.
\end{equation}
The other three canonical periods in the formula \ref{eq:PeriodsCan} can be found by solving the Picard-Fuchs operator $\mathcal{D}$ using the Frobenius method \cite{Yang}. The coefficients that occur in the prepotential $\mathcal{F}$ \ref{eq:Prepotential} have also been computed in \cite{Candelas}
\begin{equation} \label{eq:prepotentialcoeff}
Y_{111}=24,~Y_{011}=0,~Y_{001}=-2,~Y_{000}=48 \frac{\zeta(3)}{(2 \pi i)^3}.
\end{equation}
Therefore, the transformation matrix $S$ in the formula \ref{eq:smatrix} is determined up to a non-zero rational multiple by $\lambda \in \mathbb{Q}^\times$.

The authors of \cite{Candelas} have found that the smooth fiber $\mathscr{X}_{-1/7}$ of the deformation \ref{eq:aeszfamilies} over the rational point $\varphi=-1/7$ is a rank-2 attractor. More concretely, they have numerically shown that there exist two non-zero constants $C_+$ and $C_-$ such that
\begin{equation} \label{eq:AESZ34attractorEquation}
\begin{aligned}
C_+\,  S \cdot \varpi(-1/7)+\overline{C}_+ \, \overline{S} \cdot \overline{\varpi}(-1/7)&= A_+,\\
C_-\,  S \cdot \varpi(-1/7)+\overline{C}_-\, \overline{S} \cdot \overline{\varpi}(-1/7)&= A_-,
\end{aligned}
\end{equation}
where the two charge vectors $A_+$ and $A_-$ are given by \cite{Candelas}
\begin{equation} \label{eq:attractorchargesAESZ34}
A_+=(16,-60,0,5)^\top,~A_-=(0,0,2,1)^\top.
\end{equation}
Therefore the smooth fiber $\mathscr{X}_{-1/7}$ is a rank-2 attractor and the four dimensional pure Hodge structure $(H^3(X,\mathbb{Q}),F^p_{-1/7})$ at $\varphi=-1/7$ splits into the following direct sum
\begin{equation}
(H^3(X,\mathbb{Q}),F^p_{-1/7})=\mathbf{M}_{\text{B}}^{\textbf{att}}  \oplus \mathbf{M}_{\text{B}}^{\textbf{ell}},
\end{equation}
where the Hodge type of $\mathbf{M}_{\text{B}}^{\textbf{att}} $ is $(3,0)+(0,3)$, and that of $\mathbf{M}_{\text{B}}^{\textbf{ell}}$ is $(2,1)+(1,2)$. The underlying rational vector space of the direct summand $\mathbf{M}_{\text{B}}^{\textbf{att}} $ is spanned by the two charges $\beta \cdot A_+$ and $\beta \cdot A_-$ of $H^3(X,\mathbb{Z})$
\begin{equation}
\mathbf{M}_{\text{B}}^{\textbf{att}}  = \langle \beta \cdot A_+,  \beta \cdot A_- \rangle.
\end{equation}

\begin{remark}
In the next section, we will construct a natural basis for $ \mathbf{M}_{\text{B}}^\textbf{ell}$, which is also a subspace of $H^3(X,\mathbb{Q})$.
\end{remark}

\subsection{The zeta functions and \texorpdfstring{$L$}{L}-function of the rank-2 attractor  } \label{sec:aesz34Lfunction}

The authors of \cite{Candelas} have numerically computed the zeta functions of the \'etale cohomology $H^3_{\text{\'et}}(\mathscr{X}_{-1/7},\mathbb{Q}_\ell)$ for small prime numbers. At a good prime $p$, the zeta function of $H^3_{\text{\'et}}(\mathscr{X}_{-1/7},\mathbb{Q}_\ell)$ is of the form
\begin{equation}\label{eq:zetafnfactors}
(1-a_p (pT)+p(pT)^2)(1-b_p T+p^3T^2).
\end{equation}
Here, $a_p$ is the $p$-th coefficient of the $q$-expansion of a weight-2 modular form $f_2$ for the modular group $\Gamma_0(14)$, which is labeled as \textbf{14.2.a.a} in LMFDB. While $b_p$ is the $p$-th coefficient of the $q$-expansion of a weight-4 modular form $f_4$ also for the modular group $\Gamma_0(14)$, which is labeled as \textbf{14.4.a.a} in LMFDB. This property has been numerically checked by them for small prime numbers \cite{Candelas}.

The Hodge conjecture combined with the factorization of zeta functions in the formula \ref{eq:zetafnfactors} suggest that the pure motive $h^3(\mathscr{X}_{-1/7})$ splits over $\mathbb{Q}$ into the direct sum \cite{Candelas, Shamit}
\begin{equation} \label{eq:splitofAESZ34}
h^3(\mathscr{X}_{-1/7})=\mathbf{M}^{\textbf{att}} \oplus \mathbf{M}^{\textbf{ell}},
\end{equation}
where the Hodge realization of the pure motive $\mathbf{M}^{\textbf{att}}$ (resp. $\mathbf{M}^{\textbf{ell}}$) is $\mathbf{M}^\textbf{att}_{\text{B}}$ (resp. $\mathbf{M}^\textbf{ell}_{\text{B}}$). The \'etale realization of $\mathbf{M}^{\textbf{att}}$ (resp. $\mathbf{M}^{\textbf{ell}}$) is a two dimensional Galois representation whose zeta function at a good prime $p$ is $1-b_p T+p^3T^2$ (resp. $1-a_p (pT)+p(pT)^2$).

In the paper \cite{Candelas}, the special values $L(f_2,1)$, $L(f_4,1)$ and $L(f_4,2)$ have been numerically computed to a high decision
\begin{equation}
\begin{aligned}
L(f_2,1)&=0.33022365934448053902826194612283487754045234078189 \cdots ,\\
L(f_4,1)&= 0.67496319716994177129269568273091339919322842904407 \cdots, \\
L(f_4,2)&=0.91930674266912115653914356907939249680895763199044  \cdots.
\end{aligned}
\end{equation} 
There is another interesting number $v^\perp$ with numerical value
\begin{equation} \label{eq:vperp}
v^\perp= 0.37369955695472976699767292752499463211766555651682 \cdots.
\end{equation}
It has been numerically shown that the $j$-value of $\tau^\perp:=\frac{1}{2}+v^\perp i$ is rational
\begin{equation}
j(\tau^\perp)=\left(\frac{215}{28} \right)^3.
\end{equation}
The LMFDB includes only one rationally defined elliptic curve with the above $j$-invariant, which also has $f_2$ (\textbf{14.2.a.a}) as its associated eigenform under the modularity theorem of elliptic curves \cite{Candelas}. In fact, this elliptic curve is the modular curve $X_0(14)$ 
\begin{equation} \label{eq:modularcurve14}
X_0(14):y^2+xy+y=x^3+4x-6.
\end{equation}

\section{The algebraic de Rham cohomology of the rank-2 attractor } \label{sec:deRhamsplit}

In this section, we will study the algebraic de Rham cohomology of the rank-2 attractor $\mathscr{X}_{-1/7}$ and the de Rham realization of the split \ref{eq:splitofAESZ34}. More concretely, we will numerically show that the algebraic de Rham cohomology $H^3_\text{dR}(\mathscr{X}_{-1/7})$ (together with its Hodge filtration) also admits a split that `corresponds to' the split \ref{eq:splitofAESZ34}, which is induced by the two attractor equations \ref{eq:AESZ34attractorEquation}. In particular, we will construct the de Rham realizations of the pure motives $\mathbf{M}^{\textbf{att}}_{\text{B}}$ and $\mathbf{M}^{\textbf{ell}}_{\text{B}}$.

\subsection{The involution induced by complex conjugation}

First, the algebraic variety $\mathscr{X}_{-1/7}$, which is defined over $\mathbb{Q}$, also has a complex manifold structure, i.e. a Calabi-Yau manifold. The complex conjugation defines an action on the  coordinates of the points of this complex manifold, which further induces an action on its singular cohomology groups \cite{Nekovar,Yang1}. Recall that the underlying differential manifold of $\mathscr{X}_{-1/7}$ is denoted by $X$. We will denote the action of the complex conjugation on $H^3(X,\mathbb{Q})$ by
\begin{equation}
F_\infty:H^3(X,\mathbb{Q}) \rightarrow H^3(X,\mathbb{Q}).
\end{equation}  
Recall from Section \ref{sec:mirrorsymmetryattractor} that $H^3(X,\mathbb{Z})$ has an integral symplectic basis $\beta=(\beta^0,\beta^1,\alpha^0,\alpha^1)$ that satisfies the cup-product pairings in the formula \ref{eq:poincarepairing}.

The matrix of $F_\infty$ with respect to the basis $\beta$ of $H^3(X,\mathbb{Z}) $ has been computed in the paper \cite{Yang1}, which depends on the values of the coefficients $Y_{111}$, $Y_{011}$, $Y_{001}$ and $Y_{000}$ in the prepotential \ref{eq:Prepotential}. The readers are referred to the paper \cite{Yang1} for the detailed computation of $F_\infty$. Here, we plug in the values in the formula \ref{eq:prepotentialcoeff} and obtain 
\begin{equation} \label{eq:finftyMatrix}
F_\infty=
\left(
\begin{array}{cccc}
 1 & 1 & -6 & 12 \\
 0 & -1 & 12 & -24 \\
 0 & 0 & -1 & 0 \\
 0 & 0 & -1 & 1 \\
\end{array}
\right).
\end{equation}
The eigenvalues of $F_\infty$ are $1$ and $-1$, and there are two linearly independent eigenvectors associated to each eigenvalue. The two linearly independent eigenvectors of $F_\infty$ for the eigenvalue $1$ are
\begin{equation}
v^+_1=(1,0,0,0)^\top,~v^+_2=(0,-12,0,1)^\top;
\end{equation}
i.e. the subspace $H^3_+(X,\mathbb{Q})$ of $H^3(X,\mathbb{Q})$ that $F_\infty$ acts as $1$ is spanned by
\begin{equation}
\beta^0~\text{and}~-12\, \beta^1+\alpha^1.
\end{equation}
The two linearly independent eigenvectors of $F_\infty$ for the eigenvalue $-1$ are
\begin{equation}
v^-_1=(0,0,2,1)^\top,~v^-_2=(-1,2,0,0)^\top,
\end{equation}
i.e. the subspace $H^3_-(X,\mathbb{Q})$ of $H^3(X,\mathbb{Q})$ that $F_\infty$ acts as $-1$ is spanned by
\begin{equation}
2 \alpha^0+\alpha^1~\text{and}~-\beta^0+2 \beta^1.
\end{equation}

With respect to the eigenvectors of $F_\infty$, the two charge vectors $A_+$ and $A_-$ in the formula \ref{eq:attractorchargesAESZ34} can be expressed as
\begin{equation} \label{eq:chargeexpansionswrtEigen}
A_+=16v_1^++5v_2^+,~A_-=v_1^-.
\end{equation}
Recall that the two charge vectors $\beta \cdot A_+$ and $ \beta \cdot A_-$ form a basis for $\mathbf{M}^{\textbf{att}}_{\text{B}}$
\begin{equation}
\mathbf{M}_{\text{B}}^{\textbf{att}}  = \langle \beta \cdot A_+,  \beta \cdot A_- \rangle.
\end{equation}
The action of the involution $F_\infty$ on $\mathbf{M}^{\textbf{att}}_{\text{B}}$ is induced by the action of $F_\infty$ on $H^3(X,\mathbb{Q})$. From the formula \ref{eq:chargeexpansionswrtEigen}, we immediately have
\begin{equation} \label{eq:FinfinityEigenA}
F_\infty \cdot A_+=A_+,~F_\infty \cdot A_-=-A_-,
\end{equation}
which gives us the action of the involution $F_\infty$ on $\mathbf{M}^{\textbf{att}}_{\text{B}}$. We will need this action when we compute the Deligne's periods for the pure motive $\mathbf{M}^{\textbf{att}}$.

With respect to the cup-product pairing \ref{eq:poincarepairing}, the orthogonal complement of $\mathbf{M}^{\textbf{att}}_{\text{B}}$ is a two dimensional subspace spanned by
\begin{equation} \label{eq:Epm}
E_+=v_2^++6v_1^+,~E_-=v_2^- -\frac{5}{14} v_1^-.
\end{equation}
Because of the attractor equations \ref{eq:AESZ34attractorEquation}, we have
\begin{equation}
\int_X (\beta \cdot E_\pm) \smile \Omega_{-1/7}=0.
\end{equation}
The Betti realization of the pure motive $\mathbf{M}_{\text{B}}^{\textbf{ell}}$ is naturally given by
\begin{equation}
\mathbf{M}_{\text{B}}^{\textbf{ell}}  = \langle \beta \cdot E_+,  \beta \cdot E_- \rangle.
\end{equation}
From the formula \ref{eq:chargeexpansionswrtEigen}, we immediately have
\begin{equation}
F_\infty \cdot E_+=E_+,~F_\infty \cdot E_-=-E_-,
\end{equation}
which gives us the action of the involution $F_\infty$ on $\mathbf{M}^{\textbf{att}}_{\text{B}}$. We will call $\mathbf{M}^{\textbf{ell}}$ the elliptic sub-motive because of its relations to elliptic curves.

\begin{remark}
The involution $F_\infty$ \ref{eq:finftyMatrix} and its eigenvectors will play a crucial role in this paper. But we will only use them to study the arithmetic geometry of the attractor $\mathscr{X}_{-1/7}$. It is very interesting to ask whether they admit any interpretations on the K\"ahler side. Or even to ask whether do they admit any physics interpretations.
\end{remark}

\subsection{The split of the algebraic de Rham cohomology}

Let us now look at the de Rham realization $\mathbf{M}^{\textbf{att}}_{\text{dR}}$ of the pure motive $\mathbf{M}^{\textbf{att}}$ in the split \ref{eq:splitofAESZ34}. As $\mathbf{M}^{\textbf{att}}$ is a sub-motive of $h^3(\mathscr{X}_{-1/7})$, $\mathbf{M}^{\textbf{att}}_{\text{dR}}$ can be realized as a subspace of $H^3_\text{dR}(\mathscr{X}_{-1/7})$. Moreover, the Hodge filtration $F^p(\mathbf{M}^{\textbf{att}}_{\text{dR}})$ can be realized as a subspace of $F^p(H^3_\text{dR}(\mathscr{X}_{-1/7}))$. Now recall that  the Hodge filtration of the algebraic de Rham cohomology $H^3_\text{dR}(\mathscr{X}_{-1/7})$ is explicitly given by formula \ref{eq:HodgefiltrationMS}. Since the Hodge type of $\mathbf{M}^{\textbf{att}}_{\text{B}}$ is $(3,0)+(0,3)$, the Hodge filtration of $\mathbf{M}^{\textbf{att}}_{\text{dR}}$ satisfies
\begin{equation}
F^3(\mathbf{M}^{\textbf{att}}_{\text{dR}})=F^2(\mathbf{M}^{\textbf{att}}_{\text{dR}})=F^1(\mathbf{M}^{\textbf{att}}_{\text{dR}})=\langle \Omega_{-1/7} \rangle.
\end{equation}
Recall that $\langle \Omega_{-1/7} \rangle$ means the linear space spanned by $ \Omega_{-1/7} $ over $\mathbb{Q}$. 
To find the Hodge filtration of the de Rham realization $\mathbf{M}^{\textbf{att}}_{\text{dR}}$, we only need to determine $F^0(\mathbf{M}^{\textbf{att}}_{\text{dR}})$, which is a two dimensional subspace of $F^0(H^3_\text{dR}(\mathscr{X}_{-1/7}))$.

In order to construct $F^0(\mathbf{M}^{\textbf{att}}_{\text{dR}})$, we will need to find the elements of $H^3_\text{dR}(\mathscr{X}_{-1/7})$ that are orthogonal to $\mathbf{M}_{\text{B}}^{\textbf{ell}} $ under the cup-product pairing. Using numerical method, we have found that 
\begin{equation}
\int_X (\beta \cdot E_\pm) \smile \left(\Omega'''_{-1/7}-\frac{1141}{32} \Omega''_{-1/7}+\frac{15337}{64}\Omega'_{-1/7} \right)=0,
\end{equation}
hence we deduce that
\begin{equation}
F^0(\mathbf{M}^{\textbf{att}}_{\text{dR}})= \langle \Omega_{-1/7}, \Omega'''_{-1/7}-\frac{1141}{32} \Omega''_{-1/7} +\frac{15337}{64}\Omega'_{-1/7} \rangle.
\end{equation}

Let us now construct the de Rham realization $\mathbf{M}^{\textbf{ell}}_{\text{dR}}$ of the sub-motive $\mathbf{M}^{\textbf{ell}}$ in the split \ref{eq:splitofAESZ34}. First, we need to find the elements of $H^3_\text{dR}(\mathscr{X}_{-1/7})$ that are orthogonal to $\mathbf{M}_{\text{B}}^{\textbf{att}} $ under the cup-product pairing. Numerically, we have found that
\begin{equation} \label{eq:ellipticderham}
\begin{aligned}
\int_X ( \beta \cdot A_\pm ) \smile \left(\Omega'_{-1/7}-\frac{35}{8} \Omega_{-1/7} \right)&=0,\\
\int_X (\beta \cdot A_\pm) \smile \left(\Omega''_{-1/7}-\frac{735}{16} \Omega_{-1/7} \right)&=0.
\end{aligned}
\end{equation}
Since the Hodge type of $\mathbf{M}^{\textbf{ell}}_{\text{B}}$ is $(2,1)+(1,2)$, we learn that $F^3(\mathbf{M}^{\textbf{ell}}_{\text{dR}})$ must be zero. From the formula \ref{eq:ellipticderham} and the Griffiths transversality, the Hodge filtration $F^p(\mathbf{M}^{\textbf{ell}}_{\text{dR}})$ of $\mathbf{M}^{\textbf{ell}}_{\text{dR}}$ is given by
\begin{equation} \label{eq:ellipticdeRhamfiltration}
\begin{aligned}
F^2(\mathbf{M}^{\textbf{ell}}_{\text{dR}})&=\langle \Omega'_{-1/7}-\frac{35}{8} \Omega_{-1/7}\rangle, \\
F^1(\mathbf{M}^{\textbf{ell}}_{\text{dR}})&=\langle \Omega'_{-1/7}-\frac{35}{8} \Omega_{-1/7},\Omega''_{-1/7}-\frac{735}{16} \Omega_{-1/7} \rangle.
\end{aligned}
\end{equation}
This Hodge filtration $F^p(\mathbf{M}^{\textbf{ell}}_{\text{dR}})$ will be important later when we study the period of the pure motive $\mathbf{M}^{\textbf{ell}}_{\text{dR}}$.

\section{The verification of Deligne's conjecture} \label{sec:verificationDeligne}

In this section, we will numerically verify that the Tate twists of the pure motives $\mathbf{M}^{\textbf{att}}$ and $\mathbf{M}^{\textbf{ell}}$ in the split \ref{eq:splitofAESZ34} satisfy Deligne's conjecture. First, the readers are referred to the paper \cite{Yang1} for a brief review of Deligne's conjecture. 

\subsection{The attractive sub-motive} \label{sec:attractivesubmotive}

The Hodge type of the pure motive $\mathbf{M}^{\textbf{att}}$ in the split \ref{eq:splitofAESZ34} is $(3,0)+(0,3)$, hence the Tate twist $\mathbf{M}^{\textbf{att}} \otimes \mathbb{Q}(n)$ is critical if and only if $n=1,2$ \cite{DeligneL, Yang1}. Deligne's conjecture claims that the Deligne's period $c^+(\mathbf{M}^{\textbf{att}} \otimes \mathbb{Q}(n))$ is a rational multiple of the special value $L(\mathbf{M}^{\textbf{att}} \otimes \mathbb{Q}(n),s)$ when $n=1,2$.

Let us now compute the Deligne's period $c^+(\mathbf{M}^{\textbf{att}} \otimes \mathbb{Q}(n))$ for the pure motive $\mathbf{M}^{\textbf{att}} \otimes \mathbb{Q}(n)$. From the papers \cite{DeligneL, Yang1}, we have
\begin{equation}
\begin{aligned}
c^+(\mathbf{M}^{\textbf{att}} \otimes \mathbb{Q}(2))&=(2\pi i)^2\, c^+(\mathbf{M}^{\textbf{att}}),\\
c^+(\mathbf{M}^{\textbf{att}} \otimes \mathbb{Q}(1))&=(2\pi i)\, c^-(\mathbf{M}^{\textbf{att}}).
\end{aligned}
\end{equation}
Thus we only need to compute the Deligne's periods for the pure motive $\mathbf{M}^{\textbf{att}}$. The details of the computations can be found in the paper \cite{Yang1}. From the equation \ref{eq:FinfinityEigenA}, the Deligne's periods $c^\pm(\mathbf{M}^{\textbf{att}})$ for $\mathbf{M}^{\textbf{att}}$ are given by 
\begin{equation} \label{eq:attractordeligne0}
\begin{aligned}
c^+(\mathbf{M}^{\textbf{att}})&=\frac{1}{(2 \pi i)^3} \int_X (\beta \cdot A_+) \smile \Omega_{-1/7},\\
c^-(\mathbf{M}^{\textbf{att}})&=\frac{1}{(2 \pi i)^3} \int_X (\beta \cdot A_-) \smile \Omega_{-1/7}.\\
\end{aligned}
\end{equation}
To evaluate the two expressions, first recall that under the comparison isomorphism between Betti and algebraic de Rham cohomology, the threeform $\Omega_{-1/7}$ has an expansion 
\begin{equation} \label{eq:expansionOmega17}
\Omega_{-1/7}=\beta \cdot \Pi(-1/7)=\beta \cdot S \cdot \varpi(-1/7).
\end{equation}
From the cup-product pairing \ref{eq:poincarepairing}, we immediately obtain
\begin{equation} \label{eq:attractordeligne1}
\begin{aligned}
c^+(\mathbf{M}^{\textbf{att}})&=\lambda\left(-11 \varpi_0(-\frac{1}{7}) +60 \varpi_1(-\frac{1}{7}) -60 \varpi_2(-\frac{1}{7}) \right),\\
c^-(\mathbf{M}^{\textbf{att}})&=\lambda \left( \left(1-\frac{32 \zeta(3)}{(2 \pi i)^3} \right)\varpi_0(-\frac{1}{7}) + 2 \varpi_1(-\frac{1}{7}) -12 \varpi_2(-\frac{1}{7}) +  8 \varpi_3(-\frac{1}{7})\right),
\end{aligned}
\end{equation}
where $\lambda$ is the nonzero rational number in the formula \ref{eq:smatrix}. But Deligne's periods are only well-defined up to a nonzero rational multiple, so we have the freedom to simply let $c^\pm(\mathbf{M}^{\textbf{att}})$ be
\begin{equation} \label{eq:attractordeligne}
\begin{aligned}
c^+(\mathbf{M}^{\textbf{att}})&=-11 \varpi_0(-\frac{1}{7}) +60 \varpi_1(-\frac{1}{7}) -60 \varpi_2(-\frac{1}{7}) ,\\
c^-(\mathbf{M}^{\textbf{att}})&=\left(1-\frac{32 \zeta(3)}{(2 \pi i)^3} \right)\varpi_0(-\frac{1}{7}) + 2 \varpi_1(-\frac{1}{7}) -12 \varpi_2(-\frac{1}{7}) +  8 \varpi_3(-\frac{1}{7}).
\end{aligned}
\end{equation}
The numerical value of $\varpi_i(-1/7)$ can be evaluated to a high precision by numerically solving the Picard-Fuchs equation \ref{eq:aesz34PF} \cite{Candelas, Yang}. The readers are referred to the paper \cite{Yang} for more details about this method. 

From Section \ref{sec:aesz34Lfunction}, the $L$-function of $\mathbf{M}^{\textbf{att}} \otimes \mathbb{Q}(n)$ is given by \cite{DeligneL}
\begin{equation}
L(\mathbf{M}^{\textbf{att}} \otimes \mathbb{Q}(n),s)=L(f_4,s+n).
\end{equation} 
Deligne's conjecture predicts that $c^+(\mathbf{M}^{\textbf{att}} \otimes \mathbb{Q}(2))$ (resp. $c^+(\mathbf{M}^{\textbf{att}} \otimes \mathbb{Q}(1))$) is a rational multiple of $L(f_4,2)$ (resp. $L(f_4,1)$) \cite{DeligneL, Yang1}. Now we will numerically verify these two predictions. Now plug in the numerical value of $\varpi_i(-1/7)$, we immediately obtain
\begin{equation}
\begin{aligned}
c^+(\mathbf{M}^{\textbf{att}} \otimes \mathbb{Q}(2))&=-392 \, L(f_4,2), \\
c^+(\mathbf{M}^{\textbf{att}} \otimes \mathbb{Q}(1))&=14 \, L(f_4,1).
\end{aligned}
\end{equation}
Thus we have numerically shown that the critical motives $\mathbf{M}^{\textbf{att}} \otimes \mathbb{Q}(2)$ and $\mathbf{M}^{\textbf{att}} \otimes \mathbb{Q}(1)$ both satisfy Deligne's conjecture \cite{Yang1}. As Deligne's conjecture has not been proved, more examples will be important and might provide need insights into the conjecture itself.

\subsection{The elliptic sub-motive}

Let us now look at the elliptic sub-motive $\mathbf{M}^{\textbf{ell}}$, and check whether it satisfies the prediction of Deligne's conjecture. Since the Hodge type of $\mathbf{M}^{\textbf{ell}}$ is $(2,1)+(1,2)$, the Tate twist $\mathbf{M}^{\textbf{ell}} \otimes \mathbb{Q}(n)$ is critical if and only if $n=2$ \cite{DeligneL,Yang1}. From the results of Section \ref{sec:deRhamsplit}, the Deligne's periods of $\mathbf{M}^{\textbf{ell}}$ are given by \cite{DeligneL, Yang1}
\begin{equation} \label{eq:attractordeligneperiods}
\begin{aligned}
c^+(\mathbf{M}^{\textbf{ell}})&=\frac{1}{(2 \pi i)^3} \int_X (\beta \cdot E_+) \smile \left(  \Omega'_{-1/7}-\frac{35}{8} \Omega_{-1/7}\right),\\
c^-(\mathbf{M}^{\textbf{ell}})&=\frac{1}{(2 \pi i)^3} \int_X (\beta \cdot E_-) \smile \left( \Omega'_{-1/7}-\frac{35}{8} \Omega_{-1/7} \right),
\end{aligned}
\end{equation}
the numerical values of which can be computed by the similar method as in Section \ref{sec:attractivesubmotive}. While from the papers \cite{DeligneL, Yang1}, we have
\begin{equation}
\begin{aligned}
c^+(\mathbf{M}^{\textbf{ell}} \otimes \mathbb{Q}(2))&=(2\pi i)^2 \,c^+(\mathbf{M}^{\textbf{ell}}),\\
c^-(\mathbf{M}^{\textbf{ell}} \otimes \mathbb{Q}(2))&=(2\pi i)^2 \,c^-(\mathbf{M}^{\textbf{ell}}).
\end{aligned}
\end{equation}

From Section \ref{sec:aesz34Lfunction}, the $L$-function of the pure motive $\mathbf{M}^{\textbf{ell}} \otimes \mathbb{Q}(n)$ is given by
\begin{equation}
L(\mathbf{M}^{\textbf{ell}} \otimes \mathbb{Q}(n),s)=L(f_2,s+n-1).
\end{equation} 
Deligne's conjecture predicts that $c^+(\mathbf{M}^{\textbf{ell}} \otimes \mathbb{Q}(2))$ is a rational multiple of $L(f_2,1)$ \cite{DeligneL, Yang1}. Numerically, we have found that
\begin{equation} \label{eq:f2values1}
c^+(\mathbf{M}^{\textbf{ell}} \otimes \mathbb{Q}(2))=-\frac{1029}{2} \cdot L(f_2,1).
\end{equation}
Notice that as in Section \ref{sec:attractivesubmotive}, we have thrown away the rational constant $\lambda$ in the expression \ref{eq:f2values1}. This shows that the critical motive $\mathbf{M}^{\textbf{ell}} \otimes \mathbb{Q}(2)$ indeed satisfies Deligne's conjecture. 
\begin{remark}
Numerically we have also found that
\begin{equation} \label{eq:perpvalues}
c^-(\mathbf{M}^{\textbf{ell}} \otimes \mathbb{Q}(2))=-\frac{147}{8}\, \cdot \frac{L(f_2,1)\,i}{v^\perp},
\end{equation}
which provides another interesting interpretation to the number $v^\perp$ \cite{Yang1}.
\end{remark}

\section{The lattice and the complex elliptic curve} \label{sec:constructionEllipticsubmotive}

In this section, we will carefully study the elliptic sub-motive $\mathbf{M}^{\textbf{ell}}$ in the split \ref{eq:splitofAESZ34}. As the Hodge type of $\mathbf{M}^{\textbf{ell}}$ is $(2,1)+(1,2)$, so the Hodge type of the Tate twist $\mathbf{M}^{\textbf{ell}} \otimes \mathbb{Q}(1)$ is $(1,0)+(0,1)$, which is exactly the same as that of elliptic curves. In this section, we will find a two dimensional sub-lattice of $H^3(X, \mathbb{Z})$, the complex curve associated to which is the modular curve $X_0(14)$.

From Section \ref{sec:deRhamsplit}, the Betti realization $\mathbf{M}_{\text{B}}^{\textbf{ell}}$ of $\mathbf{M}^{\textbf{ell}}$ is a two dimensional subspace of $H^3(X, \mathbb{Q})$, which is spanned by $ \beta \cdot E_+$ and $ \beta \cdot E_- $. Here recall that $E_\pm$ are defined by the formula \ref{eq:Epm}. From Section \ref{sec:deRhamsplit}, the rational vector space of its Betti realization $\mathbf{M}^{\textbf{ell}}_{\text{B}} \otimes \mathbb{Q}(1)$ has a basis given by \cite{DeligneL, Yang1}
\begin{equation}
(2 \pi i) \beta \cdot E_+,~\text{and}~(2 \pi i) \beta \cdot E_-.
\end{equation}
While from the formula \ref{eq:ellipticdeRhamfiltration}, the Hodge filtration of the de Rham realization $\mathbf{M}^{\textbf{ell}}_{\text{dR}} \otimes \mathbb{Q}(1)$ is given by
\begin{equation} \label{eq:ellipticHodgefiltration}
\begin{aligned}
F^1(\mathbf{M}^{\textbf{ell}}_{\text{dR}}  \otimes \mathbb{Q}(1))&=\langle \Omega'_{-1/7}-\frac{35}{8} \Omega_{-1/7}\rangle, \\
F^0(\mathbf{M}^{\textbf{ell}}_{\text{dR}}  \otimes \mathbb{Q}(1))&=\langle \Omega'_{-1/7}-\frac{35}{8} \Omega_{-1/7},\Omega''_{-1/7}-\frac{735}{16} \Omega_{-1/7} \rangle.
\end{aligned}
\end{equation}

To construct a complex elliptic curve $\mathcal{E}$ associated to the complex elliptic sub-motive $\mathbf{M}^{\textbf{ell}} \otimes \mathbb{Q}(1)$, we will need a lattice structure for $\mathbf{M}^{\textbf{ell}}_{\text{B}} \otimes \mathbb{Q}(1)$. The integral cohomology group $H^3(X,\mathbb{Z})$ defines a lattice structure on $H^3(X,\mathbb{Q})$, which further determines a lattice structure on $\mathbf{M}^{\textbf{ell}}_{\text{B}} \otimes \mathbb{Q}(1)$. More precisely, let the two vectors $E_1$ and $E_2$ be
\begin{equation}
\begin{aligned}
E_1&=E_+=(6,-12,0,1)^\top, \\
E_2&=\frac{1}{2} \left( 5E_++14E_-\right)=(8,-16,-5,0)^\top,
\end{aligned}
\end{equation}
both of which are integral, i.e. $\beta \cdot E_i \in H^3(X,\mathbb{Z})$. A natural lattice structure on $\mathbf{M}^{\textbf{ell}}_{\text{B}} \otimes \mathbb{Q}(1)$ is generated by the two the following two elements
\begin{equation} \label{eq:latticebasis}
(2 \pi i) \beta \cdot E_1~\text{and}~(2 \pi i) \beta \cdot E_2.
\end{equation}

From the Hodge filtration of $\mathbf{M}^{\textbf{ell}}_{\text{dR}}  \otimes \mathbb{Q}(1)$  in the formula \ref{eq:ellipticHodgefiltration}, up to a nonzero constant, the nowhere-vanishing oneform of such a complex elliptic $\mathcal{E}$ must be
\begin{equation} \label{eq:holomorphiconeformE}
\Omega'_{-1/7}-\frac{35}{8} \Omega_{-1/7}.
\end{equation}
To proceed, let us compute the $j$-invariant of the elliptic curve $\mathcal{E}$ corresponding to the lattice structure of $\mathbf{M}^{\textbf{ell}}_{\text{B}} \otimes \mathbb{Q}(1)$ generated by the basis \ref{eq:latticebasis}. The period $\tau$ of $\mathcal{E}$ is given by the quotient of the pairings between the basis \ref{eq:latticebasis} and the form \ref{eq:holomorphiconeformE}, i.e.
\begin{equation}
\tau=\frac{\frac{1}{(2 \pi i)^3} \int_X (2 \pi i) \beta \cdot E_2 \smile \left(\Omega'_{-1/7}-\frac{35}{8} \Omega_{-1/7} \right)  }{\frac{1}{(2 \pi i)^3} \int_X (2 \pi i) \beta \cdot E_1 \smile \left(\Omega'_{-1/7}-\frac{35}{8} \Omega_{-1/7} \right)}.
\end{equation}
Our numerical computations have shown that $\tau$ is
\begin{equation}
\tau=\frac{5}{2}+i \cdot 0.668986610627117173003570488224 \cdots.
\end{equation}
The $j$-invariant of $\tau$ is exactly the one found in the paper \cite{Candelas}
\begin{equation}
j(\tau)=\left(\frac{215}{28} \right)^3,
\end{equation}
which is exactly that of the modular curve $X_0(14)$. 

However, $\mathcal{E}$ is only a complex elliptic curve, and $X_0(14)$ is a rational model of it. A very interesting question is whether there exists an explicit recipe that could naturally generate a rational model ($X_0(14)$) for the complex elliptic curve $\mathcal{E}$. More generally, we can ask whether physics theories can be used in the construction of such a rational model. Some progress has been made in the papers \cite{Shamit, KNY}.

\section{Conclusion and further prospects} \label{sec:conclusion}

In this paper, we have studied the connections between the rank-2 attractors and Deligne's conjecture on the special values of $L$-functions. In particular, we have focused on the example $\mathscr{X}_{-1/7}$ from the paper \cite{Candelas}. The Hodge conjecture predicts that the pure motive $h^3(\mathscr{X}_{-1/7})$ splits into the direct sum of two sub-motives $\mathbf{M}^\textbf{att}$ and $\mathbf{M}^\textbf{ell}$ over a number field $K$, where the Hodge type of $\mathbf{M}^\textbf{att}$ (resp. $\mathbf{M}^\textbf{ell}$) is $(3,0)+(0,3)$ (resp. $(2,1)+(1,2)$). While from the numerical results about the zeta functions of $h^3(\mathscr{X}_{-1/7})$ in \cite{Candelas}, $K$ should be $\mathbb{Q}$. 

In this paper, we have shown that the de Rham realization of $h^3(\mathscr{X}_{-1/7})$, i.e. the algebraic de Rham cohomology $H^3_\text{dR}(\mathscr{X}_{-1/7})$, does split accordingly over $\mathbb{Q}$, which provides further evidence to the conjectured split \ref{eq:splitofAESZ34}. Based on the numerical results of the paper \cite{Candelas}, we have computed the Deligne's periods of the critical motives $\mathbf{M}^\textbf{att} \otimes \mathbb{Q}(2)$, $\mathbf{M}^\textbf{att} \otimes \mathbb{Q}(1)$ and $\mathbf{M}^\textbf{ell} \otimes \mathbb{Q}(2)$. Then we have numerically shown that they do satisfy the predictions of Deligne's conjecture. 

We have also found a rank-2 lattice in $H^3(X, \mathbb{Z})$, which induces a lattice structure for the Betti realization $\mathbf{M}_{\text{B}}^{\textbf{ell}}$ of $\mathbf{M}^{\textbf{ell}}$. This lattice structure together with the de Rham realization $\mathbf{M}_{\text{dR}}^{\textbf{ell}}$ of $\mathbf{M}^{\textbf{ell}}$ allow us to construct a complex elliptic curve whose $j$-invariant is the same as that of $X_0(14)$. But we could not find a method to construct a rational model for this complex elliptic curve whose zeta functions agree with the results of the paper \cite{Candelas}. An open question is to ask whether physics theories can be of any help in the construction of such a rational model.

There are many interesting open questions left unanswered. For example, the computations in \ref{sec:deRhamsplit} and \ref{sec:verificationDeligne} are all on the complex side, so an interesting question is to ask whether these computations have any interpretations on the K\"ahler side. In particular, the Deligne's period $c^+(\mathbf{M}^{\textbf{att}})$ is given by the value of the period $-11 \varpi_0+60 \varpi_1-60 \varpi_2$ at $\varphi =-1/7$. So does the period $-11 \varpi_0+60 \varpi_1-60 \varpi_2$ have any interpretations on the K\"ahler side? Similarly, $c^-(\mathbf{M}^{\textbf{att}})$ is given by the value of 
\begin{equation} \label{eq:conPeriod}
\left(1-\frac{32 \zeta(3)}{(2 \pi i)^3} \right)\varpi_0 + 2 \varpi_1-12 \varpi_2+  8 \varpi_3
\end{equation}
at $\varphi = -1/7$. Does the period \ref{eq:conPeriod} admit any interpretations on the K\"ahler side?

In the paper \cite{YangK3}, we have used the method in this paper and the paper \cite{Yang1} to study the Fermat quartic 
\begin{equation}
x_0^4+x_1^4+x_2^4+x_3^4 =0,
\end{equation}
and we have shown it also satisfies the prediction of Deligne's conjecture. In the paper \cite{YangAttractor}, we have numerically shown that for $n=3,4,6,8,10$, the Fermat type Calabi-Yau $n$-fold
\begin{equation}
\sum_{i=0}^{n+1}x_i^{n+2}=0
\end{equation}
is a general rank-2 attractor. An very interesting question is to find out whether the results in this paper can be generalized to the Fermat type Calabi-Yau $n$-fold.

%%\section*{Acknowledgments}


\begin{thebibliography}{99}






\bibitem{PhilipXenia}

P. Candelas, X. C. de la Ossa, P. Green and L. Parkes. A Pair of Calabi-Yau Manifolds as an Exactly Soluble Superconformal Theory, Nuclear Physics B359 (1991) 21-74.

\bibitem{Candelas}

P. Candelas, X. de la Ossa, M. Elmi and D. van Straten. A One Parameter Family of Calabi-Yau Manifolds with Attractor Points of Rank Two. arXiv:1912.06146.

\bibitem{CoxKatz}

D. Cox and S. Katz. Mirror Symmetry and Algebraic Geometry, American Mathematical Society. 


\bibitem{DeligneL}

P. Deligne, Valeurs de fonctions $L$ et p\'eriodes d'int\'egrales. Proceedings of Symposia in Pure Mathematics 33, (1979), Part 2, 313-346.

\bibitem{Diamond}

F. Diamond and J. Shurman. A First Course in Modular Forms. Springer, GTM 228.

\bibitem{Ferrara1995} 
 
S. Ferrara, R. Kallosh and A. Strominger, N=2 extremal black holes, Phys. Rev. D. 52, R5412 (1995). DOI:10.1103/PhysRevD.52.R5412. arXiv: hep-th/9508072.
  
\bibitem{Ferrara1996} 
 
S. Ferrara and R. Kallosh, Supersymmetry and attractors, Phys. Rev. D. 54, 1514 (1996). arXiv: hep-th/9602136.

\bibitem{Yui}

F. Gouvea and N. Yui. Rigid Calabi-Yau threefolds over $\mathbb{Q}$ are modular, Expos. Math. 29 (2011), pp. 142-149.

\bibitem{MarkGross}

M. Gross, D. Huybrechts and D. Joyce. Calabi-Yau Manifolds and Related Geometries. Springer.


\bibitem{Hatcher}

A. Hatcher. Algebraic Topology, Cambridge University Press, 2008.

\bibitem{Shamit}

S. Kachru, R. Nally and W. Yang. Supersymmetric Flux Compactifications and Calabi-Yau Modularity. arXiv:2001.06022.

\bibitem{KNY}

S. Kachru, R. Nally and W. Yang. Flux Modularity, F-Theory, and Rational Models. arXiv:2010.07285.

\bibitem{KimYang}

M. Kim and W. Yang. Mirror symmetry, mixed motives and $\zeta(3)$. arXiv:1710.02344.

\bibitem{Moore}

G. Moore. Arithmetic and Attractors, arXiv:hep-th/9807087.

\bibitem{Nekovar}

J. Nekov\'a$\check{\text{r}}$. Beilinson's conjectures. Motives (Seattle, WA, 1991), 537-570, Proc. Sympos. Pure Math.,55 Part 1, Amer. Math. Soc., Providence, RI, 1994.

\bibitem{PetersSteenbrink}

C. Peters and J. Steenbrink, Mixed Hodge Structures, Springer-Verlag, 2008.


\bibitem{Yang}

W. Yang. Periods of CY $n$-folds and mixed Tate motives, a numerical study. arXiv:1908.09965.

\bibitem{Yang1}

W. Yang. Deligne's conjecture and mirror symmetry. Nuclear Physics B. 2020. DOI: \url{10.1016/j.nuclphysb.2020.115245}. arXiv: 2001.03283.

\bibitem{YangK3}

W. Yang, K3 mirror symmetry, Legendre family and Deligne's conjecture for Fermat quartic. arXiv:2004.00820.

\bibitem{YangAttractor}

W. Yang, Rank-2 attractors and Fermat type CY $n$-folds. arXiv: 2005.06722.












\end{thebibliography}
\end{document}